\documentstyle[12pt]{article}
\newcommand{\const}{\mathop{\rm const}\limits}

\newcommand{\supp}{\mathop{\rm supp}\limits}

\newcommand{\Log}{\mathop{\rm Log}\limits}

\begin{document}
\begin{center}

{\bf Wirtinger-type inequalities }\\

\vspace{3mm}

  {\bf for some rearrangement invariant spaces }\\

\vspace{3mm}

{\sc Ostrovsky E., Rogover E. and Sirota L.}\\

\normalsize

\vspace{3mm}
{\it Department of Mathematics and Statistics, Bar-Ilan University,
59200, Ramat Gan, Israel.}\\
e-mail: \ galo@list.ru \\

\vspace{3mm}

{\it Department of Mathematics and Statistics, Bar-Ilan University,
59200, Ramat Gan, Israel.}\\
e - mail: \ rogovee@gmail.com \\

\vspace{3mm}

{\it Department of Mathematics and Statistics, Bar-Ilan University,
59200, Ramat Gan, Israel.}\\
e - mail: \ sirota@zahav.net.il \\

\end{center}

\vspace{4mm}

 {\it Abstract.} In this short paper we generalize the classical inequality
 between the norms in Lebesgue spaces of the functions and its derivatives,
 which in the multidimensional case are called Sobolev's inequalities,
 on the many popular classes pairs of rearrangement invariant
(r.i.) spaces, namely, on the so-called moment rearrangement invariant spaces. \par

\vspace{3mm}


 \vspace{3mm}


\vspace{3mm}

\section{Introduction. Notations. Statement of problem.}

 Firs of all we recall "an inequality ascribed to Wirtinger (\cite{Mitrinovich1}, p. 66-68):

 $$
 \int_a^b f^2(x) dx \le \left( \frac{b-a}{2 \pi} \right)^2 \int_a^b (f^/)^2(x) dx,
 $$
 or equally

 $$
 |f|_{2,(T)} \le \frac{b-a}{2\pi} \ |f^/|_{2,(T)}. \eqno(0)
 $$
 Here $ a,b = \const, -\infty < a < b < \infty, \ T = (a,b), $ the function
 $ f(\cdot) $ has a generalized square integrable first derivative and

 $$
 f(a) = f(b), \ \int_a^b f(x) dx = 0.
 $$
 P.R. Beesack in \cite{Beesack1} obtained the following generalization of Wirtinger inequality: if $ p > 1, u^/ \in C[0,\pi/2], u(0)=0, $ then

$$
\int_0^{\pi/2} |u(x)|^p dx \le \frac{1}{p-1}
\left(\frac{p/2}{\sin(\pi/p) } \right)^p \int_0^p |u^/(x)|^p dx.
$$

  There are many generalizations of inequality (0), for example \cite{Brink1},
  \cite{Mitrinovich1}, p. 80-81:

  $$
  |f|_p  \le A(n,k) \ \Delta^{n+1/p - 1/q} |f^{(n)}|_q,  \Delta:= b-a, \
  A(n,k) < \infty,  \eqno(1)
  $$
 but in (1) the function $ f(\cdot) $ has a $ k $ fold zero at the point $ a $
 and $ (n-k) $ fold zero at the point  $ b. $ \par
 The set of all such a functions will be denoted by $ Z(n,k); 0 < k \le n: $

 $$
 Z(n,k) \stackrel{def}{=} \{f: f^{(i)}(a)=0, i=0,1,\ldots,k-1; \
 f^{(j)}(b) = 0, j=0,1,\ldots, n-k-1 \}.
 $$
 Hereafter $ n \ge 2, \ 0 \le k < n. $ \par
Evidently, the function $ f(\cdot) $ has  $ n $ times generalized derivative belonging to the space $ L_q. $ \par
 More exactly, the constants $ A(n,k) $ may be define as follows:

$$
A(n,k) = \sup_{p \in (1,\infty)} \sup_{q \in (1,\infty)}
\sup_{f \in Z(n,k), f^{(n)} \ne 0} \frac{|f|_p}{|f^{(n)}|_q} < \infty. \eqno(2)
$$

Another version of Wirtinger's inequality see, e.g. in \cite{Brink1},
\cite{Mitrinovich1}, p. 86-91: if $ f(a) = f(b)=0 $ and $ f^/ \in L_q, \ q \in (1,\infty), $ then

$$
|f|_p \le K(p,q) \ |f^/|_q, \ p \in (1,\infty), \eqno(3)
$$
where

$$
K(p,q) = \frac{q}{2} \ \frac{(1+p^*/q)^{1/p}}{(1+q/p^*)^{1/q} } \
\frac{\Gamma(1/q + 1/p^*)}{\Gamma(1/q)\Gamma(1/p^*)}, \eqno(4)
$$
$ p^* = p/(p-1)  $ and $ \Gamma(\cdot) $ denotes usually Gamma-function. \par

Note that the inequality (3) is the particular case of inequality (1) with the exact value of the constant $ A(n,k) = A(2,1). $ \par

\vspace{3mm}
In the articles \cite{Takahasi1}, \cite{Takahasi2} are considered some generalizations
of Wirtinger's inequality. In the article \cite{Giova1} was obtained the evaluated value
of the constant $ A(n,k) $ in the case of weight $ L_p - L_q $ spaces.\par

\vspace{3mm}
{\bf Our aim is generalization of Wirtinger's-type inequalities (1), (3) on some
popular classes of rearrangement invariant (r.i.) spaces, namely, on the so-called moment r.i. spaces. \par
 We intend to show also the invarianteness of offered estimations under the dilation
transform $ f \to T_{\theta}[f](x) = f(x/\theta), \ \theta = \const > 0, $
or as a minimum to show the uniform  exactness of obtained estimations at }
  $ \theta \in (0,\infty). $  \par
 The norms estimations for integral transforms, in particular, singular integral
 operators with the weight,  which are generalization of the classical Hardy-Littlewood- Weil-Rieman operators, in the Bilateral Grand Lebesgue Spaces is considered in \cite{Ostrovsky6}. \par

 \vspace{4mm}

 Hereafter $ C, C_j $ will denote any non-essential finite positive constants. As usually, for the measurable function $ f: [a,b] \to R $  we denote for sake of simplicity

 $$
|f|_p =  \left[\int_a^b |f(x)|^p \ dx \right]^{1/p}, \ 1 \le p < \infty,
 $$
  $ L_p = \{f: \ |f|_p < \infty \}; m $
 will denote usually Lebesgue measure, and we will write $ m(dx) = dx; $
  $ |f|_{\infty} \stackrel{def}{=} \sup_{x \in (a,b)} |f(x)|. $ \par
We will denote the {\it normalized } Lebesgue measure on the interval $ (a,b) $
with the length $ \Delta = b-a $ by $ m_{\Delta}: $

$$
m_{\Delta}(A) = m(A)/\Delta
$$
and will denote the correspondent $ L_p(m_{\Delta}) $ norm by $ |f|^{(\Delta)}_p: $

$$
|f|^{(\Delta)}_p = \left[ \int_a^b |f(x)|^p \ m_{\Delta}(dx)  \right]^{1/p} =
\Delta^{-1/p}|f|_p.
$$

We define also for the values $ (p_1, p_2), $ where $ 1 \le p_1 < p_2 \le \infty $

$$
L(p_1, p_2) =  \cap_{p \in (p_1, p_2)} \ L_p.
$$

\vspace{4mm}

   The Wirtinger's inequality play a very important role in the theory of approximation, theory of Sobolev's spaces, theory of function of several variables, functional analysis (imbedding theorems for Besov spaces). See,
for example,   \cite{Bennet1}, \cite{Maz'ya1}, \cite{Ostrovsky7}  etc. \par

\vspace{3mm}

  The inequality (1) may be rewritten as follows. Let $ (X, ||\cdot||X) $ be any
  rearrangement invariant (r.i.) space on the set $ T; $ denote by $ \phi(X, \delta) $ its fundamental function

  $$
  \phi(X,\delta) = \sup_{A, m(A) \le \delta} ||I(A)||X, \ I(A)= I(A,x) = I(x \in A) = 1, x \in A,
  $$
  $ \delta \in (0, \infty); \ I(A) = I(A,x) = I(x \in A) = 0, \ x \notin A.$ \par
Let us define for arbitrary r.i. space  $ (X, ||\cdot||X) $ over the set
$ (a,b) = (0,\Delta) $ the following functional:

$$
R(f; X,\Delta) \stackrel{def}{=} \frac{||f||X }{\phi(X,\Delta)},
$$
and define also for two functions
  r.i. spaces $ (X, \ ||\cdot||X ) $ and $ (Y, \ ||\cdot||Y) $
  over our set $ T = (a,b)$ with $ \Delta = b-a \in (0,\infty) $
  the so-called {\sc Wirtinger two-space functional, briefly: $ W $ functional }
  between the spaces $ X $ and $ Y $ as
  $$
  W_{n,k}(X,Y; \Delta) \stackrel{def}{=} \sup_{f \in Z(n,k), f^{(n)} \ne 0 }
 \left[  \frac{||f||X}{\phi(X,\Delta)} :
 \frac{\Delta^n \ ||f^{(n)}||Y}{\phi(Y,\Delta)} \right] =
  $$

$$
\sup_{f \in Z(n,k), f^{(n)} \ne 0 }
[R(f;X,\Delta): (\Delta^n \ R(f^{(n)}; Y,\Delta))],  \eqno(5)
$$
 or if we replace the Lebesgue measure $ m $ by the normed  measure  $ m_{\Delta} $ in the definition of the r.i. spaces $ X $ and $ Y $  and denote the correspondent norm  in the $ X,Y $ spaces over the measure $ m_{\Delta} $  by $ ||f||_{(\Delta)}X, \ ||f||_{(\Delta)}Y: $

$$
  W^{(\Delta)}_{n,k}(X,Y; \Delta) \stackrel{def}{=} \sup_{f \in Z(n,k), f^{(n)} \ne 0 } \left[\frac{||f||_{(\Delta)}X}{ \Delta^n \ ||f^{(n)}||_{(\Delta)}Y} \right]=
$$

$$
\sup_{f \in Z(n,k), f^{(n)}\ne 0}
\frac{R(f;X,\Delta)}{(\Delta^n R(f^{(n)},Y,\Delta)}. \eqno(6)
$$

  Then (1) is equivalent to the following inequalities:

  $$
  \sup_{p \in (1,\infty)} \sup_{q \in (1,\infty)}
  \sup_{ \Delta \in (0,\infty) } W^{(\Delta)}_{n,k}(L_q, L_p;\Delta) = A(n,k) < \infty, \eqno(7)
  $$

  $$
  \sup_{p \in (1,\infty)} \sup_{q \in (1,\infty)}
  \sup_{ \Delta \in (0,\infty) } W_{n,k}(L_q, L_p;\Delta) = A(n,k) < \infty. \eqno(8)
  $$

 \vspace{3mm}

{\bf Definition 1.} \par

\vspace{3mm}

{\sc By definition, the {\it pair} of r.i. spaces $ (X, \ ||\cdot||X ) $ and
  $ (Y, \ ||\cdot||Y) $ is said to be a} {\it (strong) Wirtinger's pair}, write:
  $ (X,Y) \in Wir, $ {\sc if the $ W $ functional between $ X $ and
  $ Y $  over the space $ (a,b); m $  is uniformly bounded:}
$$
\sup_{\Delta \in (0,\infty)} W_{n,k}(X,Y; \Delta) < \infty, \eqno(9)
$$

{\sc and is called a} {\it weak Wirtinger's pair}, {\sc write $ (X,Y) \in wWir, $ if for some non-trivial constant  }
$ C = C(n,k) = \const \in (0,\infty)  $

$$
\sup_{\Delta \in (0,\infty)} W_{n,k}^{(\Delta)}(X,Y; C \Delta) < \infty. \eqno(10)
$$

   {\bf Our aim is description of some pair of r.i. spaces with strong and weak
   Wirtinger properties. } \par
    Roughly speaking, we will prove that the many of popular {\it pairs} of r.i.
   spaces are strong, or at last weak Wirtinger's pairs. \par

\vspace{3mm}

  The paper is organized as follows. In the next section we recall the definition
and some properties of the so-called moment rearrangement invariant spaces, briefly,
m.r.i. spaces, which are introduced in the article \cite{Ostrovsky7} and are applied
in the theory of approximation.\par
 In the section 3 we formulate and prove the main result of this paper for m.r.i. spaces.
 In the section 4 we investigate the invariantness of obtained estimations.
 In the section 5 we will receive the Wirtinger's inequality for (generalized)
 Zygmynd spaces. \par
 The sixth section is devoted to the obtaining of the low bound for
weak Wirtinger's inequality in an arbitrary Orlicz's spaces. \par
  The last section contains some concluding remarks. \par

\vspace{4mm}

\section{Auxiliary facts. Moment rearrangement spaces.}\par

\vspace{4mm}

  Let $ (X, \ ||\cdot||X) $ be a r.i. space, where $ X $ is linear subset on the space of all measurable function $ T \to R $ over our measurable space
 $ (T,m) $ with norm $ ||\cdot||X. $ \par

\vspace{3mm}

 {\bf Definition 2.} \par

\vspace{3mm}

 {\sc We will say that the space $ X $ with the norm $ ||\cdot||X $ is {\it moment rearrangement invariant space,} briefly: m.r.i. space, or
$ X =(X, \ ||\cdot||X) \in m.r.i., $
 if there exist a real constants $ A, B; 1 \le A < B \le \infty, $ and some {\it rearrangement invariant norm } $ < \cdot > $ defined on the space of a real functions defined on the interval $ (A,B), $ non necessary to be finite on all the functions, such that

  $$
\forall f \in X \ \Rightarrow || f ||X = < \ h(\cdot) \ >, \ h(p) = |f|_p, 
\ p \in (A,B). \eqno(11)
  $$

   We will say that the space $ X $ with the norm $ ||\cdot||X $ is {\it weak moment
 rearrangement space,} briefly, w.m.r.i. space, or $ X =(X, \ ||\cdot||X) \in w.m.r.i.,$ if there exist a constants $ A, B; 1 \le A < B \le \infty, $ and some {\it functional } $ F, $ defined on the space of a real functions defined on the
interval $ (A,B), $ non necessary to be finite on all the functions, such that }

  $$
\forall f \in X \ \Rightarrow || f ||X = F( \ h(\cdot) \ ), \ h(p) = |f|_p, 
 p \in (A,B).  \eqno(12)
  $$

 We will write for considered w.m.r.i. and m.r.i. spaces $ (X, \ ||\cdot||X) $

   $$
       (A,B) \stackrel{def}{=} \supp(X),
   $$
   ("moment support"; not necessary to be uniquely defined). \par
      It is obvious that arbitrary m.r.i. space is r.i. space.\par

   There are many r.i. spaces satisfied the definition of m.r.i. or w.m.r.i
spaces: exponential Orlicz's  spaces, some Martzinkiewitz spaces, interpolation spaces (see \cite{Astashkin1},  \cite{Astashkin2}).\par

    In the article \cite{Lukomsky1} are introduced the so-called $ G(p,\alpha) $
  spaces consisted on all the measurable function $ f: T \to R $  with finite norm
   $$
   ||f||_{p,\alpha} = \left[ \int_1^{\infty} \left(\frac{|f|_x}{x^{\alpha}}
    \right)^p \ m(dx) \right]^{1/p}.
   $$

     Astashkin in \cite{Astashkin2} proved that the space $ G(p,\alpha) $
     coincides with the Lorentz $ \Lambda_p( \log^{1-p \alpha}(2/s) ) $ space.    Therefore, both this spaces are m.r.i. spaces.\par
     Another examples. Recently,, see \cite{Kozachenko1}, \cite{Fiorenza1},
     \cite{Fiorenza2}, \cite{Fiorenza3}, \cite{Iwaniec1}, \cite{Iwaniec2},
     \cite{Ostrovsky1}, \cite{Ostrovsky2}, \cite{Ostrovsky3}, \cite{Ostrovsky4},
     \cite{Ostrovsky5}, \cite{Ostrovsky6}  etc.
     appears the so-called Grand Lebesgue Spaces $ GLS = G(\psi) =
    G(\psi; A,B), \ A,B = \const, A \ge 1, A < B \le \infty, $ spaces consisting
    on all the measurable functions $ f: T \to R $ with finite norms

     $$
     ||f||G(\psi) \stackrel{def}{=} \sup_{p \in (A,B)} \left[ |f|_p /\psi(p) \right]. \eqno(13)
     $$

      Here $ \psi(\cdot) $ is some continuous positive on the {\it open} interval
    $ (A,B) $ function such that

     $$
     \inf_{p \in (A,B)} \psi(p) > 0, \ \sup_{p \in (A,B)} \psi(p) = \infty.
     $$
      It is evident that $ G(\psi; A,B) $ is m.r.i. space and
      $ \supp(G(\psi(A,B)) = (A,B). $\par

We can suppose without loss of generality
$$
\inf_{p \in (A,B)} \psi(p) = 1.
$$

  This spaces are used, for example, in the theory of probability \cite{Talagrand1}, \cite{Kozachenko1}, \cite{Ostrovsky1}; theory of Partial Differential Equations \cite{Fiorenza2}, \cite{Iwaniec2}; functional analysis \cite{Ostrovsky4}, \cite{Ostrovsky5}; theory of Fourier series \cite{Ostrovsky7}, theory of martingales \cite{Ostrovsky2} etc.\par

Note that if $ (X, ||\cdot||X) $ is m.r.i. space with the correspondent
functional  $ h(\cdot) $ and with the support $ (A,B), $ then the fundamental
function of this space has a view:

$$
\phi(\delta,X) = < f(\cdot)>, \ f(p) = \delta^{1/p}, \ p \in (A,B).
$$
For instance, the fundamental function for the Grand Lebesgue Space $ G(\psi) $
with  the support $ (A,B) $ may be calculated by the formula

$$
\phi(\delta,G(\psi)) = \sup_{p \in (A,B)} \frac{\delta^{1/p}}{\psi(p)}.
$$
 The detail investigation of fundamental functions for Grand Lebesgue Spaces,
with  consideration of many examples, see, e.g. in \cite{Ostrovsky8}. \par

\vspace{3mm}

   Let us consider now the (generalized) Zygmund's spaces $ L_p \ (\Log)^r L, $ which may 
   be defined as an Orlicz's spaces over some subset of the space $ R^l $  with non-empty interior  and with $ N - $ Orlicz's function of a view

    $$
    \Phi(u) = |u|^p \ \log^r( C + |u|), \ p \ge 1, \ r \ne 0.
    $$

    It is known \cite{Ostrovsky7} that: \\
    {\bf 1.} All the spaces $ L_p \ (\Log)^r L $ over real line with measure
  $ m $  with condition $ \ r \ne 0 $ are not m.r.i. spaces.\\
     {\bf 2.} If $ r $ is positive and integer, then the spaces $ L_p (\Log)^r L $   are w.m.r.i. space.\\
 {\bf 3.} There exists an r.i. space without the w.m.r.i. property.\par

\vspace{3mm}

\section{Main result. Wirtinger's inequality for the pairs of m.r.i. spaces.}

\vspace{3mm}

{\bf Theorem 3.1.}\par
{\it Let $ (X, ||\cdot||X) $ be any m.r.i. space
 relatively the auxiliary norm $ < \cdot >, $ and let $ (Y,||\cdot||Y ) $ be
another m.r.i. space over at the same set $ (T,m)$  relatively the
second auxiliary norm $ << \cdot >>. $ \par

  Then the pair of m.r.i. spaces $ (X, ||\cdot||X) $ and $ (Y, ||\cdot||Y) $
is the (strong) Wirtinger's pair uniformly in } $ \Delta, \
\Delta \in (0,\infty): $

$$
\sup_{\Delta > 0} W_{n,k}(X,Y) \le A(n,k) < \infty. \eqno(14)
$$

{\bf Proof} is very simple.  Let $ f $ be arbitrary function from the set
$ Z(n,k): f \in Z(n,k) $ and let $ f^{(n)} \ne 0. $

 It follows from the Brink's inequality (1) that

$$
|f|_p  \ \Delta^{1/q} \le A(n,k) \ \Delta^n
|f^{(n)}|_q \ \Delta^{1/p}. \eqno(15)
$$
We get tacking the norm $ << \cdot >> $ from both sides of inequality (15):

$$
|f|_p \ \phi(Y, \Delta) \le A(n,k) \ \Delta^n ||f^{(n)}||Y \cdot \Delta^{1/p}.
\eqno(16)
$$
We obtain now tacking the norm $ < \cdot > $ from both sides of inequality (16):

$$
||f||X \ \phi(Y,\Delta) \le A(n,k) \ \Delta^n \ ||f^{(n)}||Y \cdot
\phi(X,\Delta),
$$
which is equivalent to the assertion of the considered theorem.\par
Note as an example that for the Grand Lebesgue Spaces $ G(\psi) $ and $ G(\nu) $
the proposition (14) may be rewritten as follows. Let us denote

$$
V_{\Delta}(f; G(\nu), G(\psi))= \left[ \frac{||f||G(\nu)}{\phi(G(\nu),\Delta)}:
\Delta^n \frac{||f^{(n)}||G(\psi)}{\phi(G(\psi), \Delta)} \right];
$$
then

$$
\sup_{\Delta > 0} \sup_{f \in Z(n,k), f^{(n)} \ne 0} V_{\Delta}(f; G(\nu),
G(\psi)) \le A(n,k) < \infty.  \eqno(17).
$$

\vspace{3mm}

\section{Low bounds for Grand Lebesgue Spaces.}

\vspace{3mm}

We investigate in this section the exactness of inequality (12), or in other words the asymptotical invariableness  under the dilation operators $ T_{\Delta}.$ \par
 Note that

 $$
 (T_{\Delta}f)^{(n)} = \Delta^n T_{\Delta}f^{(n)}.
 $$

  Let $ g: [0,1] \to R $ be any function from the set $ Z(n,k) $ such that
  $ g^{(n)} \ne 0. $  We continue this function at the values
 $ x \ge 1 $ by zero: $ x > 1 \ \Rightarrow g(x) = 0. $ \par

 Let us denote

$$
V_0(\psi,\nu)= \inf_{\Delta \in (0,\infty)} \sup_{g \in Z(n,k), g^{(n)} \ne 0}
\left[ \frac{||T_{\Delta}g||G(\nu)}{\phi(G(\nu),\Delta)}:
\frac{||T_{\Delta}g^{(n)}||G(\psi)}{\phi(G(\psi),\Delta)}  \right].
$$

{\bf Theorem 4.1.}

$$
V_0(\psi,\nu) \ge \frac{k^k \ (n-k)^{n-k}}{n! \ (n+1)^{n+1}}. \eqno(18)
$$
{\bf Proof.} We note first of all that if a function $ g: R_+ \to R $ is such that
 for some positive finite constants $ C_1 $ and $ C_2 $

 $$
 C_1 \le \inf_{p \in [0,1]} |g|_p \le \sup_{p \in [1,\infty]} |g|_p \le C_2,
 $$
then

$$
R(T_{\Delta}g; G(\psi), \Delta) \in [C_1, C_2].
$$

Further, let us choose

$$
\alpha = \frac{k}{n+1}, \ \beta = \frac{n-k}{n+1},
$$
and
$$
g(x):= g_{n,k}(x) = x^k \ (1-x)^{n-k}, \ x \in [0,1]; \ g_{n,k}(x) = 0, x > 1.
$$

We conclude after simple calculations that when $ x \in [\alpha,\beta], $ then

$$
g_{n,k}(x) \ge  \frac{k^k \ (n-k)^{n-k} }{ (n+1)^n },
$$
and

$$
\max_{x \in [0,1]} |g_{n,k}(x)| = \frac{k^k \ (n-k)^{n-k}}{n^n};
$$
therefore

$$
\frac{k^k \ (n-k)^{n-k}}{(n+1)^{n+1}} \le |g_{n,k}|_p \le
\frac{k^k \ (n-k)^{n-k}}{n^n}.
$$

Substituting into the expression for the value  $ V_0(\psi,\nu), $ we obtain the
inequality

$$
V_0(\psi,\nu) \ge \frac{k^k \ (n-k)^{n-k}}{n! \ (n+1)^{n+1}},
$$
Q.E.D.\par
 Note as a corollary that we obtain the following {\it low} bound for the constant
 $ A(n,k):$
 $$
 A(n,k) \ge \frac{k^k \ (n-k)^{n-k}}{n! \ (n+1)^{n+1}}. \eqno(19)
 $$

\vspace{3mm}

\section{The case of (generalized) Zygmund spaces. Other method.}

\vspace{3mm}

 Recall that the (generalized) Zygmund space

$$
 X = L_q \ (\Log)^{\gamma} \ L
$$
or correspondingly $ Y = L_p \ (\Log)^{-\beta} \ L $ is defined as an Orlicz's space
with the Orlicz's function of a view:

$$
 \Phi(u) = \Phi(q, \gamma; u) = |u|^q \ [\log(C(q,\gamma) + u)]^{\gamma}, \eqno(20)
$$
where $ C(q, \gamma) $ is sufficiently great positive constant.\par
 We assume in this section that  $ p > 1, \ (q < \infty), \ \beta,\gamma > 0. $ \par

 Note that the fundamental functions for these spaces are as $ \Delta
 \in (0, \infty): $

$$
\phi(L_q \ (\Log)^{\gamma} \ L), \Delta) \asymp
 \Delta^{-1/q} \ (1 + |\log \Delta|)^{\gamma/q}. \eqno(21)
$$
 Let $ Y $ be another Zygmund's space: $ Y = L_p \ (\Log \ L)^{- \beta}. $
 We will formulate and prove now the
Wirtinger's inequality for Zygmund spaces, but {\it only in the cases}
$  \gamma, \beta \ge 0.$ \par
 Let us denote

$$
 L_{q,+} =  \cup_{\epsilon \in (0,1)} L_{q + \epsilon}
$$
and correspondingly

$$
L_{p,-} = \cup_{ \delta \in (0, 0.5(p+1))}L_{p - \delta };
$$
we define also for the measurable function $ f: [0,1] \to R $ and $ f(x) = 0,
 x > 1 $ the following quotient (Wirtinger's functional):

$$
W^o(\Delta; p,q; n,k) =  \sup_{f \in Z(n,k), f \in L_{q,+}}\sup_{f \in L_{p,-}}
\left[\frac{||T_{\Delta}f||X}{\phi(X,\Delta)}:
\frac{\Delta^n ||T_{\Delta} f^{(n)}||Y}{\phi(Y,\Delta)} \right]. \eqno(22)
$$

{\bf Theorem 5.1.} {\it Let } $ \gamma \ge 0, \beta \ge 0. $
{\it We assert that Wirtinger's functional for considered spaces is uniformly over
the variable } $ \Delta $ {\it is bounded:}

$$
 \sup_{\Delta \in (0,\infty) } W^o(\Delta; p,q; n,k) = C(p,q;n,k) < \infty.
$$

{\bf Proof.} Since the cases $ \gamma = 0 $ or $ \beta = 0 $ are simple, we investigate further only the possibility $ \gamma > 0, \ \beta > 0. $ \par

 It is proved the articles \cite{Grand1}, \cite{Nessel1} that for
$ g \in L_{q,+}  $ and for arbitrary values $ r \in (q, q +1) $

$$
||g||L_q \ (\Log \ L)^{\gamma} \le C \left[ \frac{r}{r - q} \right]^{\gamma/r} \
|g|_r \eqno(23)
$$
and analogously may be proved the 'inverse' inequality: for arbitrary \\
$ s \in (0.5(1+ p),p) $

$$
||g||L_p \ (\Log)^{-\beta} L \ge C \left[ \frac{s}{p - s} \right]^{- \beta/s} \
|g|_s. \eqno(24)
$$

 We have for $ g \in L_{p,+} $ and for the values $ r \in (q, q + 1) $

$$
R(g; X,\Delta) \stackrel{def}{=}
\frac{||g||X}{\phi(X,\Delta)} \le C \cdot R_1,
$$
where
$$
 R_1 = R_1(g; X,\Delta) \stackrel{def}{=}
\frac{||g||X}{\Delta^{1/q} [1+| \log \Delta|]^{\gamma/q}} \le
$$

$$
C \frac{(r/(r-q))^{\gamma/r} \ |g|_r}{\Delta^{1/q} \ [1+|\log \Delta|]^{\gamma/q}}
$$
and we find analogously for the values $ s \in (0.5(1+p), p) $

$$
 R(g; Y,\Delta) \stackrel{def}{=}
\frac{\Delta^n ||g^{(n)}||Y}{\phi(Y,\Delta)} \ge C \cdot R_2,
$$
where
$$
 R_2 = R_2(g^{(n)}; Y,\Delta) \stackrel{def}{=}
\frac{\Delta^n ||g^{(n)}||Y}{\Delta^{1/p} [1+|\log \Delta|]^{- \beta/p}} \ge
$$

$$
C \frac{(s/(p-s))^{-\beta/s} \ |\Delta^n g^{(n)}|_s}{\Delta^{1/p} \ [1+|\log \Delta|]^{-\beta/p}}.
$$

 The assertion of theorem 5.1 may be obtained after the dividing the estimation for
 $ R_1 $ over the estimation for $ R_2, $ using the Brink's inequality (1) for the
 estimation of the quotient $ |g|_r/\Delta^n |g^{(n)}|_s $ and after the
 minimizing over $ (s,r). $ \par

  More simple, we can choose in order to prove theorem 5.1 in the expression
 for $ R_1/R_2 $ for all sufficiently greatest values $ |\log \Delta| $

$$
r := r_0 = q + \frac{\gamma}{q \ [1 + |\log \Delta|]}, \
s := s_0 = p - \frac{\beta}{p \ [1 + | \log \Delta| ]}.
$$

\vspace{3mm}

\section{Low bound for arbitrary Orlicz's spaces. }

\vspace{3mm}

Let us consider in this section a two arbitrary Orlicz's spaces
$ L(\Phi) $ and $ L(\Phi_1) $ over the set $ (a,b) = (0,\Delta)  $ with the
correspondent Orlicz's functions $ \Phi = \Phi(u) $ and $ \Phi_1 = \Phi_1(u). $ \par
 Let also $ g = g(x), \ x \in [0,1] $ be some function from the set $ Z(n,k) $
such that $ g^{(n)} \ne 0. $ Denote $ \overline{W}_{n,k} = $
$$
\overline{W}_{n,k}(\Phi,\Phi_1) =
\inf_{\Delta \in (0,\infty)} \sup_{g \in Z(n,k), g^{(n)} \ne 0}
\left[\frac{||T_{\Delta}g||L(\Phi)}{\phi(L(\Phi),\Delta/(n+1))}:
\frac{||T_{\Delta} g^{(n)}||L(\Phi_1)}{\phi(L(\Phi_1),\Delta)} \right].
$$

{\bf Theorem 6.1.}

$$
\overline{W}_{n,k} \ge \frac{k^k \ (n-k)^{n-k} }{n! \ (n+1)^n}. \eqno(25)
$$
{\bf Proof.} Recall that the norm of the measurable function $ h: (0,\Delta) \to R $ in the $ L(\Phi) $ space may be introduced, for instance, by the  formula

$$
||h||L(\Phi)= \inf_{v > 0} v^{-1} \left[1 + \int_0^{\Delta} \Phi(v h(x)) dx  \right].
 \eqno(26)
$$

 For example, if $ h = T_{\Delta}g, \ g: (0,1) \to R,$ then

$$
||T_{\Delta}g||L(\Phi) = \inf_{v > 0} v^{-1}
\left[1 + \Delta \int_0^1 \Phi(v g(y)) dy  \right].
$$
Let us choose as before

$$
g(x) = g_{n,k}(x)= x^k (1-x)^{n-k}, \ x \in [0,1],
$$
then $ g_{n,k}(\cdot) \in Z(n,k) $ and $ |g^{(n)}_{n,k}| = n!. $ \par
Therefore,

$$
||T_{\Delta} g^{(n)}_{n,k}(\cdot)||L(\Phi_1) = \inf_{v > 0} v^{-1}
[1+ \Delta \Phi_1(v \ n!)] =
$$

$$
n! \ \inf_{v > 0} v^{-1} [1+ \Delta \Phi_1(v)].
$$

 Note that the fundamental function for the $ L(\Phi) $ space has a view:

$$
\phi(L(\Phi), \delta) = \inf_{v> 0} v^{-1} \left[ 1+\delta \Phi(v) \right], \eqno(27)
$$
following

$$
R(g_{n,k}; L(\Phi_1), \Delta) = n!.
$$

Furthermore, let us choose as before

$$
\alpha = \frac{k}{n+1}, \ \beta = \frac{n-k}{n+1},
$$
we conclude that when $ x \in [\alpha,\beta], $ then

$$
g_{n,k}(x) \ge  \frac{k^k \ (n-k)^{n-k} }{ (n+1)^n },
$$
therefore

$$
||g_{n,k}(\cdot)||L(\Phi) \ge \inf_{v > 0} v^{-1} \left[ 1 + \Delta
\int_{\alpha}^{\beta} \Phi(v g_{n,k}(x) dx) \right] \ge
$$

$$
\inf_{v > 0} v^{-1} \left[ 1 + \Delta \int_{\alpha}^{\beta}
\Phi(v \min_{x \in [\alpha,\beta]}g_{n,k}(x)) \ dx \right] \ge
$$

$$
\inf_{v > 0} v^{-1} \left[ 1 + (\Delta /(n+1))
 \Phi(v k^k \ (n-k)^{n-k}/(n+1)^n) \right]=
$$

$$
\frac{k^k \ (n-k)^{n-k}}{(n+1)^n} \cdot \inf_{v > 0}v^{-1}
\left[ 1 + \Delta/(n+1) \Phi(v) \right]=
$$

$$
\frac{k^k \ (n-k)^{n-k}}{(n+1)^n} \phi(L(\Phi),\Delta/(n+1))=
\frac{k^k \ (n-k)^{n-k}}{(n+1)^n}.
\eqno(28)
$$
Dividing the last estimation on the $ \phi(L(\Phi),\Delta/(n+1)), $
 we obtain the assertion (25) of theorem 6.1. \par

\vspace{3mm}

\section{ Concluding remarks. }

\vspace{3mm}

We consider in this section some slight generalization of Wirtinger's-Brink's
inequalities  (3)-(4) on the Grand Lebesgue Spaces, with at the same notations,
for instance $ K(p,q) \ $ (4). \par
 We suppose for definiteness $ a = 0, \ b = 1, $ so that $ \Delta = 1. $ \par

Let $ f: [0,1] \to R $ be some function such that for some
 $ \psi \in \Psi(A_1, B_1) \Rightarrow  df/dx (\cdot) \in G(\psi; A_1, B_1). $ \par

Note that the function
 $ \psi(\cdot) $ may be "constructive"  introduced by means of equality

$$
\psi(q):= |df/dx(\cdot)|_q,
$$
(the so-called natural choice), if of course
 $ \psi(q) < \infty, \ q \in (A_1, B_1). $ \par
 Let us define the function $ \nu(p) $ by the following way:

$$
\nu(p) = \inf_{q \in (A_1, B_1)} [K(p,q) \psi(q)], \eqno(29)
$$
 denote
$$
(A_2, B_2) = \{p: \ \nu(p) < \infty \}
$$
and suppose $ 1 \le A_2 < B_2 \le \infty. $ \par

{\bf Theorem 7.1.} {\it Let } $ f(0) = f(1) = 0. $ {\it We assert that }

$$
||f||G(\nu) \le ||f^/||G(\psi). \eqno(30)
$$
{\bf Proof.} Let $ f(0) = f(1) = 0 $ and $ f^/ \in G(\psi); $ then

$$
|f^/|_q \le ||f^/||G(\psi) \cdot \psi(q).
$$

We get using the Brink's inequality:

$$
|f|_p \le K(p,q) \psi(q) \ ||f^/||G(\psi),
$$
and thus

$$
|f|_p \le \inf_{q \in (A_1, B_1)} [K(p,q) \ \psi(q)] \ ||f^/||G(\psi) =
\nu(p) \ ||f^/||G(\psi), \eqno(31)
$$
which is equivalent to the assertion of theorem 7.1  by virtue of definition of the
norm in $ G(\nu) $ space.\par

\end{document}